\documentclass[11pt]{article}
\begin{document}
\tolerance=10000


\hyphenpenalty=2000
\hyphenation{visco-elastic visco-elasticity}
\setcounter{page}{1}
\thispagestyle{empty}




\font\note=cmr10 at 10 truept  


\def\pni{\par \noindent}
\def\vsh{\smallskip}
\def\s{\smallskip}
\def\vs{\medskip}
\def\vvs{\bigskip}
\def\vvvs{\bigskip\medskip} 
\def\vsp{\par}
\def\vsn{\vsh\pni}
\def\cen{\centerline}
\def\ra{\item{a)\ }} \def\rb{\item{b)\ }}   \def\rc{\item{c)\ }}
\def\eg{{\it e.g.}\ } \def\ie{{\it i.e.}\ }


\def\sg{\hbox{sign}\,}
\def\sgn{\hbox{sign}\,}
\def\sign{\hbox{sign}\,}
\def\e{{\rm e}}
\def\exp{{\rm exp}}
\def\ds{\displaystyle}
\def\dis{\displaystyle}
\def\q{\quad}    \def\qq{\qquad}
\def\lan{\langle}\def\ran{\rangle}
\def\l{\left} \def\r{\right}
\def\lra{\Longleftrightarrow}
\def\d{\partial}
\def\dr{\partial r}  \def\dt{\partial t}
\def\dx{\partial x}   \def\dy{\partial y}  \def\dz{\partial z}
\def\rec#1{{1\over{#1}}}
\def\zr{z^{-1}}



\def\hatt{\widehat}
\def\epsilons{{\widetilde \epsilon(s)}}
\def\sigmas{{\widetilde \sigma (s)}}
\def\fs{{\widetilde f(s)}}
\def\Js{{\widetilde J(s)}}
\def\Gs{{\widetilde G(s)}}
\def\Fs{{\wiidetilde F(s)}}
 \def\Ls{{\widetilde L(s)}}
\def\L{{\cal L}} 
\def\F{{\cal F}} 


\def\NN{{\bf N}}
\def\RR{{\bf R}}
\def\CC{{\bf C}}
\def\ZZ{{\bf Z}} 


\def\I{{\cal I}}  
\def\D{{\cal D}}  


\def\Gc{{\cal {G}}_c}   \def\Gcs{\barr{\Gc}} 
\def\Gs{{\cal {G}}_s}   \def\Gss{\barr{\Gs}} 
\def\Gck{\hatt{\Gc}} 
\def\args{(x/ \sqrt{D})\, s^{1/2}}
\def\argsa{(x/ \sqrt{D})\, s^{\beta}}
\def\arg{ x^2/ (4\,D\, t)}


\def\erf{\hbox{erf}}  \def\erfc{\hbox{erfc}}


\def\uks{{\widehat{\widetilde {u}}} (\kappa,s)}

\def\psikappa{\psi_\alpha^\theta(\kappa)}






\cen{FRACALMO PRE-PRINT   {\bf ww.fracalmo.org}}
\vsh
\cen{Published in {\bf Vietnam Journal of Mathematics}
  Vol. 32 (2004), SI pp. 53-64}
\vsh
\cen{Springer Verlag:  ISSN 0866-7179}
\vsh
\hrule
   \vskip 0.50truecm


\font\title=cmbx12 scaled\magstep2
\font\bfs=cmbx12 scaled\magstep1

\begin{center}

{\title A fractional generalization}
\vs

{\title of the Poisson processes}

\vvs

 {Francesco MAINARDI}$^{(1)}$,
{Rudolf GORENFLO}$^{(2)}$,
{ Enrico SCALAS} $^{(3)}$


\vs

$\null^{(1)}$
 Dipartimento di Fisica, Universit\`a di Bologna and INFN, \\
Via Irnerio 46, I-40126 Bologna, Italy \\
E-mail: {\tt mainardi@bo.infn.it} $\;$ URL: {\tt www.fracalmo.org}
\\ [0.25 truecm]

$\null^{(2)}$
 First Mathematical Institute,
 Free   University of Berlin, \\
 Arnimallee  3, D-14195 Berlin, Germany \\
E-mail: {\tt gorenflo@mi.fu-berlin.de}
\\ [0.25 truecm]

$\null^{(3)}$
Dipartimento di Scienze e Tecnologie Avanzate,
Universit\`a del \\ Piemonte Orientale,
via Cavour 84,  I-15100 Alessandria, Italy\\
E-mail: {\tt scalas@unipmn.it}
 \def\date#1{\gdef\@date{#1}} \def\@date{\today}

\end{center}



\cen{\bf Abstract} 

\vskip 0.1truecm
\noindent
It is our intention to provide via fractional calculus
a  generalization of the pure and compound Poisson processes,
which are known to play a fundamental role in renewal theory,
without and with reward, respectively.
We first recall the basic renewal theory including its
fundamental concepts like waiting time between events, the survival
probability, the counting function.
If the waiting time is exponentially distributed we  have a
Poisson process, which is Markovian.
However, other waiting time distributions are also relevant in
applications, in particular such ones with a fat tail caused by a power
law decay of its density. In this context we analyze a non-Markovian
renewal process with a waiting time distribution described by the
Mittag-Leffler function.
This distribution, containing the exponential
as  particular case,  is shown to play a fundamental role
in the infinite thinning procedure of  a generic renewal process
governed by a power-asymptotic waiting time.
We then consider the renewal theory with reward that implies a
random walk subordinated to a  renewal process.

\vskip 6pt


\vskip 6pt


\vfill\eject
\section{\bfs Essentials of  renewal theory}

The concept of  {\it renewal process} has been developed as a
stochastic model for describing the class of counting processes
for which the times between successive  events are
independent  identically distributed  ($iid$)
non-negative random variables, obeying a given probability law.
These times are  referred to as waiting times
or inter-arrival times.
For more details see \eg the classical
treatises by
Khintchine \cite{Khintchine QUEUEING60},
Cox \cite{Cox RENEWAL67},
Gnedenko \& Kovalenko \cite{GnedenkoKovalenko QUEUEING68},
Feller \cite{Feller 71},
and the  recent book by Ross \cite{Ross PROBMOD97}.

For a renewal process having waiting times $T_1,T_2, \dots$, let
$$t_0 = 0\,, \q   t_k= \sum_{j=1}^k T_j\,, \q k \ge 1\,. \eqno(1.1)
    $$
That is $t_1 =T_1$ is the time of the first renewal,
$t_2 = T_1 +T_2$ is the time of the second renewal and so on.
In general $t_k$ denotes the $k$th renewal.

The process is specified if we know the probability law for the waiting
times. In this respect we introduce the
{\it probability density function} ($pdf$)
$\phi(t)$
 and  the   (cumulative) distribution function $\Phi(t)$
so defined:
$$ \phi (t) := \frac{d}{dt} \Phi(t) \,,\q
   \Phi(t) := P \l( T \le t\r) = \int_0^t \phi (t')\, dt'\,.
\eqno(1.2)$$
 When the nonnegative random variable  represents
 the lifetime of technical systems, it is common to refer
to $\Phi(t)$ as to the {\it failure probability}
and to
$$ \Psi(t) := P \l(T > t\r) = \int_t^\infty \phi (t')\, dt' = 1-\Phi(t)\,,
\eqno(1.3)$$
as to the {\it survival probability}, because $\Phi(t)$ and $\Psi(t)$ are
the respective probabilities that the system does or does not fail
in $(0, T]$. 
A relevant quantity is the {\it counting function} $N(t)$
defined as
$$ N(t):= \hbox{max} \l\{k | t_k \le t, \;k = 0, 1, 2, \dots\r\}\,,
\eqno(1.4)$$
that represents the effective number of events before or at instant $t$.
In particular we have $\Psi(t) = P \l(N(t) =0\r)\,.$
Continuing in the general theory we set
$F_1(t) = \Phi(t)$, $f_1(t) = \phi (t)$, and in general
$$ F_k(t) :=  P\l(t_k = T_1+ \dots +T_k \le t \r)\,,
 \; f_k(t) = \frac{d}{dt} F_k(t)\,, \; k\ge 1\,, \eqno(1.5)$$
thus  $F_k(t)$ represents the probability that the sum of the first $k$
waiting times is less or equal  $t$
and $f_k(t)$ its density. Then, for any fixed $k\ge 1$  the normalization
condition  for $F_k(t)$ is fulfilled
because
$$         \lim_{t \to \infty} F_k(t) =
  P\l(t_k = T_1+ \dots +T_k < \infty \r)= 1 \,.
                 \eqno(1.6)$$
In fact,    the sum of $k$ random variables each of which is
    finite with probability 1 is finite with probability 1 itself.
By setting for consistency $F_0(t) \equiv 1$ and $f_0(t) = \delta(t)$,
the Dirac delta  function\footnote{
We find it convenient to recall the {\it formal representation}
 of  this generalized function in $\RR^+\,,$
$$\delta(t) := \frac{t^{-1}}{\Gamma(0)}\,, \q t\ge 0\,.$$},
we also note that for $k \ge 0$ we have
$$ P\l(N(t) =k\r) := P\l(t_k \le t \,,\,t_{k+1} > t\r)
=  \int_0^t f_k(t')\, \Psi (t-t')\, dt'\,.
\eqno(1.7)$$

We now find it convenient to introduce the simplified $\, *\,$ notation
for the Laplace convolution between two causal well-behaved
(generalized) functions $f(t)$ and $g(t)$
$$   \int_0^t   f(t')\, g(t-t')\, dt' = \l(f \,*\, g\r) (t) =
     \l(g \,* \, f\r) (t)  = \int_0^t   f(t-t')\, g(t')\, dt'\,. $$
Being  $f_k(t)$ the $pdf$ of the sum of the  $k$
$iid$ random variables
$T_1,  \dots, T_k$
with $pdf$ $\phi (t)\,, $ we easily recognize that
$f_k(t)$ turns out to be the $k$-fold convolution of $\phi(t)$
with itself,
$$ f_k(t) =  \l(\phi^{*k}\r) (t)\,, \eqno(1.8)$$
so Eq. (1.7)  simply reads:
$$    P\l(N(t) =k\r) = \l(\phi^{*k} \,*\, \Psi\r)(t)\,. \eqno(1.9)$$
Because of the presence of Laplace convolutions
a renewal process is suited for the Laplace transform method.
Throughout this paper we will denote
by   $\widetilde f(s)$
the Laplace transform
of a sufficiently well-behaved (generalized) function
$f(t)$  according to
$$   {\L} \l\{ f(t);s\r\}=   \widetilde f(s)
 = \int_0^{+\infty} \e^{\ds \, -st}\, f(t)\, dt\,,
\q s > s_0\,,
$$
and for $\delta (t)$
consistently we will have
$  \widetilde \delta (s) \equiv 1\,. $
{\it} Note that for our purposes we agree to take $s$ real.
We recognize  that
(1.9) reads in the Laplace domain
$$    \L \{P\l( N(t)=k\r); s\}
= \l[{\widetilde \phi (s)} \r]^k \,   \widetilde \Psi (s)
             \,,\eqno(1.10)$$
where, using (1.3),
$$ \widetilde \Psi (s)  =
\frac{ 1- \widetilde \phi (s)} {s}\,.
\eqno(1.11)$$

\section{{\bfs The Poisson process as a renewal process}}

The most celebrated   
renewal process is the Poisson process
characterized by a waiting time $pdf$ of exponential type,
$$  \phi (t) = \lambda \, \e^{-\lambda t}\,,\q \lambda >0\,, \q t\ge 0\,.
\eqno(2.1)$$
The process has {\it no memory}.
Its  moments turn out to be
 $$ \langle T\rangle = \frac{1}{\lambda }\,, \q
   \langle T^2 \rangle = \frac{1}{\lambda^2 }\,,
\q \dots \,,  \q \langle T^n \rangle = \frac{1}{\lambda^n }\,, \q \dots \,,
\eqno(2.2)$$
and the {\it survival probability} is
  $$ \Psi(t) := P\l(T>t\r) = \e^{-\lambda t}\,, \q t \ge 0 \,.
\eqno(2.3)$$
We know that the probability  that $k$ events occur in the
interval of length $t $ is
$$ P\l( N(t) = k\r) =   \frac{(\lambda t)^k}{k!} \, \e^{-\lambda t}
\,, \q t \ge 0\,, \q k = 0,1, 2, \dots\;. \eqno(2.4)$$
The probability distribution related to the sum of $k$ $iid$
exponential random variables  is known to be
the so-called {\it Erlang distribution}  (of order $k$).
The corresponding density (the {\it Erlang} $pdf$) is thus
$$ f_k(t) = \lambda \,  \frac{(\lambda t)^{k-1}}{(k-1)!}
    \, \e^{-\lambda t}\,,\q t \ge 0 \,,
\q k =1,2, \dots \,, \eqno(2.5)$$
so that the  Erlang distribution function of order $k$ turns out
to be
$$ F_k(t) = \int _0^t f_k(t')\, dt' =
     1 -  \sum_{n=0}^{k-1}
      \frac{(\lambda t)^n}{n!} \, \e^{-\lambda t}  =
 \sum_{n=k}^{\infty}
      \frac{(\lambda t)^n}{n!} \, \e^{-\lambda t}\,,\q t \ge 0
\,. \eqno(2.6)$$
In the limiting case $k=0$ we recover
$f_0(t) = \delta(t),\; F_0(t) \equiv 1,\; t\ge 0$.

The results  (2.4)-(2.6) can easily obtained
by using the technique of the Laplace transform sketched in
the previous section
noting that for the Poisson process we have:
$$ \widetilde\phi(s) =
\frac{\lambda } {\lambda +s}\,,
 \q  \widetilde\Psi(s) =  \frac  {1 } {\lambda +s}\,,   \eqno(2.7)$$
and for the Erlang distribution:
$$  \widetilde f_k(s) =  [\widetilde\phi(s)]^k =
    \frac  {\lambda^k } {(\lambda +s)^k}\,,
\q    \widetilde F_k(s) =  \frac{[\widetilde\phi(s)]^k}{s} =
    \frac  {\lambda^k } {s (\lambda +s)^k}\,.
\eqno(2.8)$$

We also recall that the survival probability
for the Poisson renewal process
obeys the ordinary differential equation  (of relaxation type)
$$     \frac{d}{dt} \Psi(t) = -\lambda \Psi(t)\,, \q t \ge 0\,; \q
\Psi(0^+) =1\,. \eqno(2.9)$$

\section{{\bfs The renewal  process of Mittag-Leffler type} }

A "fractional" generalization of the   Poisson renewal process
is simply obtained by generalizing the differential equation  (2.9)
replacing there the first derivative with the integro-differential operator
$\, _tD_*^\beta$ that  is interpreted  as
the fractional derivative
of order $\beta $ in  Caputo's sense, see Appendix.
We write, taking for simplicity $\lambda =1$,
$$       \, _tD_*^\beta \, \Psi(t) =
- \Psi(t)\,, \q t >0\,,\q 0<\beta \le 1\,; \q \Psi(0^+) =1\,. \eqno(3.1)$$
We also allow
the limiting case $\beta =1$ where all the results of the previous
section (with $\lambda =1$) are expected to be recovered.

For our purpose we need to recall the Mittag-Leffler function
as the natural "fractional" generalization
of the exponential function, that characterizes the Poisson process.
The Mittag-Leffler function of parameter  $\beta\,$
is defined in the complex plane by the power series
$$ E_\beta (z) :=
    \sum_{n=0}^{\infty}\,
   {z^{n}\over\Gamma(\beta\,n+1)}\,, \q \beta >0\,, \q z \in \CC\,.
 \eqno  (3.2)$$
It turns out  to be an entire function of order $\beta $
which reduces for $\beta=1$ to $\exp (z)\,.$
For detailed information on the Mittag-Leffler-type functions
and their Laplace transforms the reader  may  consult \eg
\cite{Erdelyi HTF,GorMai CISM97,Podlubny BOOK99}.

The solution  of Eq. (3.1) is known to be, see \eg
\cite{CaputoMainardi RNC71,GorMai CISM97,Mainardi CHAOS96},
$$\Psi(t) =  E_\beta (-t^\beta)\,, \q t \ge 0\,, \q 0<\beta \le 1\,,
\eqno (3.3)$$ so
$$ \phi(t) :=    -   \frac{d} {dt} \Psi(t) =
            -   \frac{d}{dt}  E_\beta (-t^\beta) \,, \q t \ge 0
 \,, \q 0<\beta \le 1\,.\eqno (3.4)$$
Then, the corresponding Laplace transforms read
 $$ \widetilde \Psi(s) =
 {s^{\beta-1} \over 1+ s^\beta}\,,   \q
 \widetilde \phi(s)= {1 \over 1+  s^{\beta}}\,,\q
    0<\beta \le 1\,.\eqno(3.5) $$
Hereafter, we find it convenient to summarize
the  most relevant features  of the functions $\Psi(t)$ and $\phi(t)$
when  $0< \beta <1\,.$
We begin to quote their series expansions for $t \to 0^+ $
and  asymptotics for  $t\to \infty $,
$$ \Psi(t)
   = {\ds \sum_{n=0}^{\infty}}\,
  (-1)^n {\ds {t^{\beta n}\over\Gamma(\beta\,n+1)}}
 \,\sim \,  {\ds {\sin \,(\beta \pi)\over \pi}}
  \,{\ds  {\Gamma(\beta)\over t^\beta}}\,,
     \eqno(3.6) $$
and
$$ \phi(t)
= {\ds {1\over t^{1-\beta}}}\, {\ds \sum_{n=0}^{\infty}}\,
  (-1)^n {\ds {t^{\beta n}\over\Gamma(\beta\,n+\beta )}}
 \, \sim \,  {\ds {\sin \,(\beta \pi)\over \pi}}
  \,{\ds  {\Gamma(\beta+1)\over t^{\beta+1}}}\,.
     \eqno(3.7) $$
In contrast 
to the Poissonian case  $\beta=1$,
in the case  $0<\beta <1$ for large $t$
the functions $\Psi(t)$ and $\phi(t)$
no longer decay   exponentially  but algebraically.
As a consequence of the power-law asymptotics
the process turns be no longer Markovian but of long-memory type.
However, we recognize that for $0<\beta <1$ both  functions
 $\Psi(t)$, $\phi(t)$
keep   the "completely monotonic" character of the Poissonian case.
Complete monotonicity of  the   functions
 $\Psi(t)$ and $\phi(t)$  means
$$ (-1)^n {d^n\over dt^n}\, \Psi  (t) \ge 0\,,  \q
   (-1)^n {d^n\over dt^n}\, \phi  (t) \ge 0\,,
\q n=0,1,2,\dots   \,, \q t \ge 0\,,     \eqno(3.8)$$
or equivalently, their representability as real Laplace transforms
of non-negative   generalized
functions (or measures), see \eg \cite{GorMai  CISM97}. 

For the generalizations
of Eqs (2.4) and (2.5)-(2.6),  characteristic
of the Poisson and Erlang distributions respectively,
we must point out the   Laplace transform pair
$$ \L\{ t^ {\beta \,k}\, E_\beta ^{(k)}
  (-t^\beta ) ;s\} =
        \frac{ k!\, s^{\beta -1}}{(1+s^\beta )^{k+1}}
\,, \q \beta >0 \,, \q k = 0,1, 2, \dots \,,
\eqno(3.9)$$
with $ {\ds E_\beta ^{(k)}(z) := \frac{d^k}{dz^k}  E_\beta(z)}\,, $
that can be deduced  from the book by Podlubny,
see (1.80) in  \cite{Podlubny BOOK99}.
Then, by using the Laplace transform pairs  (3.5) and
Eqs (3.3), (3.4), (3.9)
in Eqs  (1.8) and (1.9),
 we have  the {\it generalized Poisson distribution},
$$ P\l( N(t) = k\r) =   \frac{ t^{ k\, \beta}}{k!} \,
  E_\beta^{(k)} (- t^\beta)
\,, \q k = 0, 1, 2, \dots \eqno(3.10)$$
and the {\it generalized Erlang} $pdf$   (of order $k \ge 1$),
$$ f_k(t) = \beta \,  \frac{ t^{k\beta-1}}{(k-1)!}
    \, E_\beta^{(k)} (- t^\beta)   
\,. \eqno(3.11)$$
The {\it generalized  Erlang distribution function} turns out
to be
$$ F_k(t) = \int _0^t f_k(t')\, dt' =
     1 -  \sum_{n=0}^{k-1}
      \frac{ t^{n \beta}}{n!} \, E_\beta^{(n)} (- t^\beta)  =
 \sum_{n=k}^{\infty}
      \frac{t^{n\beta}}{n!} \, E_\beta^{(n)} (- t^\beta)
\,. \eqno(3.12)$$


\section{{\bfs The Mittag-Leffler distribution as limit for thinned
 renewal processes}}

The procedure of  thinning (or rarefaction) for a generic renewal
process  (characterized by a generic random sequence of waiting times
$\{T_k\}$)
has been considered and investigated by Gnedenko and Kovalenko
\cite{GnedenkoKovalenko QUEUEING68}.
It  means that for each positive
index  $k$  a decision is made: the event is deleted with probability $p$
or it is maintained with probability $q=1-p$, with
 $0<q<1$.
For this thinned or rarefied renewal process
we  shall hereafter   revisit and complement the results available
in   \cite{GnedenkoKovalenko QUEUEING68}.
 We begin to rescale the time variable  $t$
replacing it by  $t/r$, with a parameter $r$
on which we will dispose later.
Denoting, like in (1.5),  by $F_k(t)$
the probability distribution function of the sum of
$k$  waiting times and by $f_k(t)$ its density,
we have recursively, in view of (1.8),
$$ f_1(t) = \phi (t) \,,\q
      f_k(t) = \int_0^t  f_{k-1}(t-t') \, \phi (t')\,dt'
  = \l(\phi^{*k}\r)(t)\,,\q  
k \ge 2\,. \eqno(4.1)$$
Let us denote by $(T_{q,r} f)(t) $
 the waiting time density 
in the thinned and rescaled process from one event to the next.
Observing that after a maintained event the next one of the
original process is kept with probability $q$ but dropped in favour
of the second next with probability $p\,q$  and, generally,
$n-1$ events are dropped in favour of the  $n$-th next with probability
$p^{n-1} \, q$ ,
we arrive at the formula
$$ (T_{q,r} f)(t) = \sum_{n=1}^\infty q\, p^{n-1}\, f_n(t/r)/r\,.
\eqno(4.2)$$

Let $\widetilde f_n(s) = \int_0^\infty \e^{-st}\, f_n(t)\, dt $
   be the Laplace transform of $f_n(t)$.
 Recalling $f_1(t)= \phi(t) $
we set $ \widetilde f_1(s)= \widetilde \phi (s)$.
Then $f_n(t/r)/r$   has the transform
 $\widetilde f_n(rs)= \l(\widetilde \phi (rs)\r)^n$,
and we obtain (in view of $p=1-q$)
the formula
$$ (T_{q,r} \widetilde \phi)(s)
= \sum_{n=1}^\infty q\, p^{n-1}\,\l(\widetilde \phi (rs)\r)^n
 = \frac{ q\, \widetilde \phi (rs)}{ 1 - (1-q)\, \widetilde \phi (rs)}
 \,,\eqno(4.3)$$
from which by Laplace inversion we can, in principle, construct
the transformed process.

  We now imagine stronger and  stronger rarefaction (infinite thinning)
by considering a scale of processes with  the parameters $r=\delta $
and $q= \epsilon $  tending
to zero under a scaling relation $\epsilon =\epsilon (\delta )$
  yet to be specified.
Gnedenko and Kovalenko have, among other things, shown
that if the condition
$$\widetilde \phi(s) =1-a(s) \, s^\beta  + o\l(a(s)\, s^\beta \r)\,,
 \q \hbox{for} \q s \to 0^+\,,\eqno(4.4)$$
where $a(s)$ is a    slowly varying function for $s\to 0$
\footnote{{\bf Definition:} We call a (measurable) positive function $a(y)$,
defined in a right neighbourhood of zero, {\it slowly varying at zero} if
$a(cy)/a(y) \to 1$ with $y \to 0$ for every $c>0$.
We call a (measurable) positive function $b(y)$,
defined in a  neighbourhood of infinity, {\it slowly varying at infinity}
if
$b(cy)/a(y) \to 1$ with $y \to \infty$ for every $c>0$.
Examples: $(\log y)^{\gamma}$ with $\gamma \in \RR\,$ and
$\,\exp \,\l({\log y}/{\log\, \log y}\r)$.}
, is satisfied, then we have
with $\epsilon =\epsilon (\delta )= a(\delta)\, \delta^\beta$
for every {\it fixed} $s>0$   the limit relation
$$
\widetilde \phi _0(s) :=
\lim_{\delta \to 0} \frac{ \epsilon (\delta )\, \widetilde \phi (\delta s)}
   {1 - (1-\epsilon )\, \widetilde \phi (\delta s)}=
   \frac{1}{ 1 + s^\beta }\,,\q 
 \q 0<\beta \le 1\,. \eqno(4.5)$$
This condition is met with
$a(s) = \lambda \, M(1/s)$ if the waiting time $T$ obeys a power
law with index $\beta $, in the sense of {\it Master Lemma 2}
by Gorenflo and Abdel-Rehim \cite{GAR Vietnam03}.
The function $M(y)$ is the same as in {\it Master Lemma 2},
so it varies slowly at infinity,
whence $M(1/s)$ varies slowly at zero.
The proof of (4.5) is by straightforward calculation.
Observe the slow variation property of $a(s)$ and note that
terms small of higher order become negligible in the limit.
By the continuity theorem for Laplace transforms,
see Feller \cite{Feller 71}, we now recognize $\phi _0(t)$
as the limiting density,  which we identify,
in view of (3.2)-(3.5),
$$
    \phi_0(t)= - \frac{d}{dt} E_\beta (- t^\beta )\,.
\eqno(4.6)$$
So the limiting waiting time density is
the so-called Mittag-Leffler density, that in the special case $\beta =1$
reduces to the well-known exponential density. 
\vfill\eject

\noindent
It should be noted that  Gnedenko and Kovalenko
in the sixties
failed to recognize  $\widetilde\phi _0(s)$ as Laplace transform
of a Mittag-Leffler type function\footnote
{Although the Mittag-Leffler function
was introduced by the Swedish mathematician G. Mittag-Leffler
in the first years of the twentieth century,
it lived for long time  in  isolation as
{\it Cinderella}. 
The  term  {\it Cinderella function}  was used in the fifties
by the Italian mathematician
F.G. Tricomi for the incomplete gamma function.
In recent years the Mittag-Leffler function is gaining
more and more popularity in view of the increasing applications
of the fractional calculus and is classified as 33E12
in the  Mathematics Subject Classification 2000.}.

\vvs

\section{{\bfs  Renewal processes with reward:
the fractional master equation and its solution}}

The renewal process can be accompanied  by reward that means
that at every renewal instant a space-like variable makes a  random jump
from its previous position to a new point in "space".
"Space" is here  meant in a very general sense.
In the insurance business, e.g., the renewal points
are instants where the company receives a payment or must
give away money to some claim of a customer,
so space is money.
In such process occurring in time and in space,
also referred to as {\it compound renewal process},
the probability distribution of jump widths is as relevant
as that of waiting times.

Let us denote by $X_n$  the jumps occurring at instants $t_n\,,$
$\, n = 1,2,3,\dots\,$.
Let us assume that $X_n$  are $iid$ (real, not necessarily positive)
 random variables
with  probability density $ w(x)$,  independent of  the
{\it waiting time} density $\phi (t)$.
In a physical context the $X_n$s represent
the jumps of a diffusing particle (the walker),
and the resulting random walk model is known
as {\it continuous time random walk} (abbreviated as CTRW)
in that the  waiting time is assumed to be
a {\it continuous}  random variable
\footnote{The name CTRW became
popular in physics after that in the 1960s
 Montroll, Weiss and Scher (just to cite  the pioneers)
published a celebrated series
of papers on random walks to model
diffusion processes on lattices, see \eg
\cite{Weiss 94} and references therein.
CTRWs are rather good and general phenomenological models for diffusion,
including anomalous diffusion, provided that the resting time of the
 walker is much greater than the time it takes to make a jump.
In fact, in the formalism, jumps are instantaneous.
In more recent times, CTRWs were applied back to economics and
finance by  Hilfer \cite{Hilfer 84}, by the authors of
the present paper with M. Raberto
\cite{Scalas 00,Mainardi 00,Gorenflo 01,Raberto 02},
and, later, by Weiss and co-workers \cite{Masoliver 03a}.
However,  it should be noted that
the idea of combining a stochastic process for waiting times between
two consecutive events and another stochastic process which associates
a reward or a claim to each event dates back at least to the first half
 of the twentieth century with the so-called
 Cram\'er--Lundberg model for insurance risk, see for a review
\cite{Embrechts 01}. In a probabilistic framework, we now find it
more appropriate to refer to all these processes as to
{\it compound renewal processes}.
}.
\vfill\eject
\noindent
The position $x$ of the walker at time $t$ is
$$
x(t) = x(0) + \sum_{k=1}^{N(t)} X_k\,.
\eqno(5.1)$$
Let us now denote by $p(x,t)$   the probability density
 of finding the random walker
at the position $x$ at time instant $t\,. $  We assume
the initial condition $p(x,0) = \delta (x)\,, $  meaning that the
walker is initially at the origin,  $x(0)=0\,. $
We look for the evolution equation for $p(x,t) $
 of the compound renewal process.
Based upon the previous probabilistic arguments
we arrive at
$$   p(x,t) =  \delta (x)\, \Psi(t) +
   \int_0^t   \phi (t-t') \, \l[
 \int_{-\infty}^{+\infty}  w(x-x')\, p(x',t')\, dx'\r]\,dt'
 \,, \eqno(5.2) $$
called the {\it integral equation of the CTRW}.
In fact,  from Eq. (5.2) we recognize the role of the
{survival probability}
$\Psi(t)$ and of the densities $\phi (t)\,,\,  w(x)\,.$
The first term in the RHS of (5.2) expresses the persistence
(whose strength decreases with   increasing time)
of the initial position $x=0$.
The second term (a space-time convolution) gives
the contribution to $p(x,t)$ from the walker sitting
in point $x' \in \RR$ at instant $t' < t$ jumping to
point $x$ just at instant $t\,,$  after stopping (or waiting) time
$t-t'\,. $

The integral equation (5.2) can be solved
by using the machinery of the transforms of Laplace
and Fourier.
Having introduced the notation for the Laplace transform
in sec. 1, we now quote
our notation for the Fourier transform,
$ \F\{f(x); \kappa\} = \widehat{f}(\kappa)=\int_{-\infty}^{+\infty}
\e^{i \kappa x} \,f(x) \,dx $ ($\kappa \in \RR$),
 and for the corresponding Fourier  convolution between
(generalized) functions
$$ \l(f_{1}\,*\, f_{2}\r)(x)= \int_{-\infty}^{+\infty}
f_{1}(x')\,f_{2}(x-x')\, dx'\,.$$
Then, applying the transforms of Fourier and Laplace in succession to
(5.2) and using the well-known
operational rules, we
arrive at 
the famous Montroll-Weiss equation, see \cite{MontrollWeiss 65},
$$
 \widehat{\widetilde{p}} (\kappa,s) = 
\frac{\widetilde\Psi(s)}{1- \tilde{\phi}(s)\,\widehat{w}(\kappa)}
\,.
\eqno(5.3)
$$
As pointed out in \cite{GAR Vietnam03},
this equation can alternatively be derived from the Cox formula, see
    \cite{Cox RENEWAL67}  chapter 8 formula (4),
describing the process as    subordination of a random
    walk to a renewal process.
\vfill\eject

By inverting the transforms one can, in
    principle, find the evolution $p(x,t)$  of the sojourn density for
    time $t$    running from zero to infinity.
In fact, recalling that $|\widehat w(\kappa)| < 1$ and
$|\widetilde\phi (s)| < 1$,
if $\kappa \not= 0$ and
$s \not= 0$, Eq. (5.3) becomes
$$
\widetilde{\widehat p}(\kappa, s) =
\widetilde \Psi(s)\, \sum_{k=0}^{\infty}
[\widetilde \phi (s) \, \widehat w(\kappa)]^k \,;
\eqno(5.4) $$
this gives, inverting the Fourier and
the Laplace transforms and taking into account
Eqs. (1.9)-(1.10),
$$
p(x,t) = \sum_{k=0}^{\infty} P(N(t)=k)\, w_k (x)\,,
\eqno(5.5)$$
where $w_k(x) = \l(w ^{*k}\r)(x),\,$
in particular $w_0(x) = \delta(x), \; w_1(x) = w(x).$

A special case of the integral  equation (5.2) is obtained for
the {\it compound Poisson process} where $\phi (t) = \e^ {-t}$
(as in (2.1) with $\lambda =1$ for simplicity).
Then, the  corresponding equation reduces after some
manipulations, that best are carried out in the Laplace-Fourier domain,
  to the {\it Kolmogorov-Feller equation}:
$$
  \frac{\d }{\d  t}\,p(x,t)= -p(x,t)+\int_{-\infty}^{+\infty}
w(x-x')\, p(x',t)\, dx' \, ,
\eqno(5.6)
$$
which is the {\it master equation of the compound Poisson process}.
In this case, in view of Eqs (2.4) and (5.5) the solution reads
$$p(x,t) = \sum_{k=0}^{\infty} \frac{ t^k}{k!} \,
{\e}^{- t} \,w_k (x)\,.
\eqno(5.7)
$$
When the survival probability is the Mittag-Leffler function
introduced in    (3.3),  the master equation
                  for the corresponding
fractional version of the compound process
can be shown to be
   $$
  \null_tD_*^\beta \,p(x,t)= -p(x,t)+\int_{-\infty}^{+\infty}
w(x-x')\, p(x',t)\, dx' \,, \q 0<\beta <1\,,
\eqno(5.8)
$$
where  $\,_tD_*^\beta$ denotes the time fractional
derivative   of order $\beta $  in the Caputo sense.
For a (detailed) derivation of Eq (5.8 ) we refer  to the paper
by Mainardi et al. \cite{Mainardi 00}, in which the results
have been obtained by an  approach independent
from that adopted in a previous paper
by Hilfer and Anton \cite{HilferAnton 95}.
In this case, in view of Eqs (3.10) and (5.5), the solution
of the {\it fractional master equation} (5.8) reads:
$$ p(x,t) = \sum_{k=0}^{\infty} \frac{ t^{\beta k}}{k!} \,
   E_\beta^{(k)} (- t^\beta)  \,w_k (x)\,.
\eqno(5.9)
$$
In \cite{Gorenflo 01}  we have, under a power law regime for
the jumps,  investigated for Eq. (5.8)
the so-called {\it diffusive or hydrodynamic limit},
obtained by making smaller
all jumps by a positive factor $h$ and
accelerating  the process  by a large factor properly related
to $h$, then letting  $h$ tend to zero.
In this limit the master equation (5.8) reduces
to  a {\it space-time fractional diffusion equation}.
This is also the topic
 of the recent paper by Scalas et al. \cite{Scalas FRACTALS03} 
and, in a  more general framework, of the paper
by Gorenflo and Abdel-Rehim \cite{GAR Vietnam03}.
\vvs

\noindent
{\bfs{Conclusions}}

\vskip 0.2truecm
\noindent
We have provided a  {\it fractional generalization}
of the Poisson renewal processes
by replacing  the first time derivative
in the relaxation equation of the survival probability
by a fractional derivative of order $\beta$
($0 <\beta \le 1$). Consequently, we have obtained for
$0<\beta<1$  non-Markovian renewal processes
where, essentially, the exponential probability 
densities, typical for the Poisson processes,
are replaced by functions of Mittag-Leffler type,
that decay in a  power law manner with an exponent related to $\beta$.

The  distributions  obtained by considering the
sum of $k$ $iid$   random variables distributed according
to the Mittag-Leffler  law provide the
fractional generalization of the corresponding
Erlang distributions.     
Furthermore, the Mittag-Leffler probability distribution
is shown to be the limiting distribution for the thinning procedure
of a generic renewal process with waiting time density
of power law character.

Then, our theory  has been  applied to renewal processes
{\it with reward}, so can be considered as
the fractional generalization of the compound Poisson processes.
In such processes, occurring in time and in space,
also the probability distribution of the jump widths is relevant.
The stochastic evolution of the space variable in time
is modelled by an    integro-differential equation
(the master equation) which,  by containing a time fractional derivative,
can be considered as the fractional generalization of the
classical Kolmogorov-Feller equation of the  compound Poisson process.
For this master equation we have provided the analytical solution
in terms of iterated derivatives of a Mittag-Leffler function.
\vfill\eject

\noindent
{\bfs{Appendix: The Caputo fractional derivative}}

\vskip 0.25 truecm
\noindent
The     {\it Caputo} fractional   derivative     
provides a fractional generalization
of the first derivative  through the following
rule in the Laplace transform domain,
$$ {\cal L} \l\{ _tD_*^\beta  \,f(t) ;s\r\} =
      s^\beta \,  \widetilde f(s)
   -  s^{\beta  -1}\, f(0^+) \,,
  \q 0<\beta  \le 1 \,,\q s>0\,,  \eqno(A.1)$$
hence  turns out to be defined as,
see \eg \cite{CaputoMainardi RNC71,GorMai CISM97},
$$
    _tD_*^\beta \,f(t) :=
\cases{
    {\ds \rec{\Gamma(1-\beta )}}\,{\ds\int_0^t
 {\ds {f^{(1)}(\tau)\over (t-\tau )^{\beta }}\, d\tau
 }} \,,
  & $\; 0<\beta  <1\,, $\cr\cr
     {\ds {d\over dt}} f(t)\,,
    & $\; \beta  =1\,. $\cr\cr }
   \eqno(A.2) $$
It can alternatively be written in the form
$$
\qq \qq \qq     _tD_*^\beta \,f(t) =
 {\ds \rec{\Gamma(1-\beta )}}\, {\ds {d \over dt}}
{\ds\int_0^t}
 {\ds {f(\tau)\over (t-\tau )^{\beta }} \, d\tau} \, - \,
{ t^{-\beta }\over \Gamma(1-\beta )} \, f(0^+)
\qq \qq \qq \eqno(A.3)$$
$$  =
    {\ds \rec{\Gamma(1-\beta )}}\, {\ds {d \over dt}}
{\ds\int_0^t}
 {\ds {f(\tau) -f(0^+) \over (t-\tau )^{\beta }} \, d\tau}
 \,,\q 0<\beta <1
 \,. $$
The Caputo derivative has been indexed with $*$
in order to distinguish it from the classical
Riemann-Liouville fractional derivative
$_tD^\beta$,
the first term at the R.H.S. of the first equality
 in (A.3). As it can be noted from the last equality in (A.3),
the Caputo derivative provides a sort of regularization at $t=0\, $
of the Riemann-Liouville derivative; for more details see
\cite{GorMai CISM97}.



\end{document}